\documentclass[12pt]{amsart}%
\usepackage{amsmath}
\usepackage{amsfonts}
\usepackage{amssymb}
\usepackage{graphicx}%
\providecommand{\U}[1]{\protect\rule{.1in}{.1in}}
\newtheorem{theorem}{Theorem}
\newtheorem{acknowledgement}{Acknowledgement}

\newtheorem{axiom}{Axiom}

\newtheorem{corollary}[theorem]{Corollary}

\newtheorem{definition}[theorem]{Definition}

\newtheorem{proposition}[theorem]{Proposition}

\newtheorem{rmk}{Remark}
\newenvironment{Rmk}{\begin{rmk}\em}{\end{rmk}}
\newtheorem{prf}{Proof}

\newenvironment{Prf}{\begin{prf}\em}{\qed\end{prf}}

\newcommand{\NOT}[1]{}
\newcommand{\pa}{\par\medskip}

\title[(Positive) Totally Ordered Noncommutative Monoids]{\textbf{(Positive) Totally Ordered Noncommutative Monoids -- How Noncommutative Can They Be?}}
\author{Eliahu Levy}
\address{Department of Mathematics\\
Technion -- Israel Institute of Technology\\
Technion City, Haifa 3200003, Israel}
\email{(eliahu@math.technion.ac.il)}
\date{May 2020}

\begin{document}

\begin{abstract}
Commutative totally ordered monoids abound, number systems for example. When the monoid is not assumed commutative,
one may be hard pressed to find an example. One suggested by Professor Orr Shalit are the countable ordinals with addition.

In this note we attempt an introductory investigation of totally (also partially) ordered monoids, not assumed commutative (still writing them additively), and taking them as positive, i.e.\ every element is greater than the unit element. That, in the usual commutative cases, allows the ordering to be defined via the algebraic structure, namely, as divisibility (in our additive sense): $a\le b$ defined as $\exists\,c\,\,(b=a+c)$. The noncommutative case offers several ways to generalize that.

First we try to follow the divisibility definition (on the right or on the left). Then, alternatively, we insist on the ordering being compatible with the operation both on the left and on the right, but strict inequality may not carry over -- again refer to the ordinals example. We try to see what axiom(s) such requirements impose on the monoid structure, and some facts are established.

Focusing especially on the totally ordered case, one finds that necessarily the noncommutativity is somewhat limited. One may partly emulate here the commutative case, speaking about infinitely grater vs.\ Archimedean to each other elements, and in the Archimedean case even emulate Euclid's Elements' theory of `ratios' -- all that imposing some partial commutativity.
\end{abstract}

\maketitle

\section{Introduction}\label{sc:int}
One encounters many examples of totally ordered \textit{commutative} monoids $M$ (which we shall write additively). They have a unit element (denoted by $0$), and have the set $M_+$ of the nonnegative elements, which is a submonoid, and the ordering may be defined via $M_+$: $b\ge a\equiv\exists c\in M_+\left(b=a+c\right)$. Further restricting ourselves to $M_+$, i.e., assuming $M=M_+$, we refer to the monoid as \textit{positive}. Then
\begin{equation}\label{eq:def0}
b\ge a\Leftrightarrow\exists c(b=a+c).
\end{equation}
I.e., speaking multiplicatively, $b\ge a$ means $b$ \textit{is divisible} by $a$.\pa

And this is assumed to be a total ordering: transitive (as (\ref{eq:def0}) implies, using the associativity of the monoid) but also antisymmetric and total: for $a\ne b$, one and only one of $a>b$ or $b>a$ holds. Of the many examples one may cite the positive integers (or rationals, reals)\pa

But what if the monoids is noncommutative? then even in (\ref{eq:def0}) there may be a difference between requiring
divisibility in the left or on the right (or maybe both?). Also, one would be hard pressed to find an indeed
noncommutative example.\pa

An example, as noted by Professor Orr Shalit, are the ordinals (say the countable ordinals, to make them a set and not a proper class), with respect to \textit{addition}. Here indeed addition is associative, but not commutative: $\omega+1\ne\omega=1+\omega$, and the ordering is characterized by ({\ref{eq:def0}}) \textit{in that order: $b=a+c$}. Indeed, by the theory of ordinals (and well-ordered sets), $b\ge a$ means that $a$ is an initial segment of $b$, and $c$ will be what $b$ has after that initial segment.\pa

Note that here we have compatibility of the ordering: $a\le b\Rightarrow c'+a\le c'+b$, which of course follows from the characterization ({\ref{eq:def0}}) (and associativity), but also compatibility on the right: $a\le b\Rightarrow a+c\le b+c$, but here strict inequality would not be carried over: $0<1$ but $0+\omega=\omega=1+\omega$.\pa

In this note it is attempted to shed some introductory light on these issues (with the ordering partial or total).\pa

In the first section we try to follow the divisibility `paradigm' in the noncommutative case, assuming cancelation (both on the left and on the right), compare \cite{MS-MS} Ch.\ 4. We show by example that one may still have then both noncommutativity and a divisibility ordering which is total.\pa

Next, alternatively, one insists on the ordering being compatible with the operation both on the left and on the right, embarking on a somewhat methodical investigation about what axiom(s) would be imposed by our requirements on the monoid structure, proving some facts. We focus then on the possibilities for idempotents and an element $a$ absorbing an element $b$ on the left or on the right (i.e.\ $a+b=a$ or $b+a=a$) in positive partially ordered monoids -- one does not find them when there is cancelation, but in the general case they may play some part.\pa

Finally the totally ordered case is considered. There one finds that necessarily the noncommutativity is somewhat limited. One may partly emulate here the commutative case, speaking about infinitely greater vs.\ Archimedean to each other elements, and in the Archimedean case even emulate Euclid's Elements' theory of `ratios' -- all that imposing some partial commutativity.

\subsection{Some Notation}
$M$ is a monoid (written additively, though not assumed commutative), i.e., a set with an associative operation $+$ having a unit element $0$. We endow it with an ordering relation $\le$ (whose inverse relation we, as usual, denote by $\ge$). We consider the following postulates that may (or may not) hold:
\begin{eqnarray*}
&&\textrm{(O) Transitive: If }a\le b\textrm{ and }b\le c\textrm{ then }a\le c,\quad\forall a,b,c\in M.\\
&&\textrm{(P) Positive: }a\ge0,\quad \forall a\in M.\\
&&\textrm{(C) Compatible: If }a\le b\textrm{ then }a+c\le b+c\textrm{ and }c+a\le c+b,\quad\forall a,b,c\in M.\\
&&\textrm{(T) Total: For any }a,b\in M,\textrm{ either }a\le b\textrm{ or }b\le a\textrm{ (allowing both)}.\\
&&\textrm{(A) Antisymmetric: For }a,b\in M,\textrm{ if both }a\le b\textrm{ and }b\le a\textrm{ then }a=b.
\end{eqnarray*}
\begin{definition}
A monoid endowed with an ordering relation satisfying (O), (P), (C) as above is called a \textbf{positive preordered monoid}. If (A) is also satisfied, a \textbf{positive partially ordered monoid}, and if both (A) and (T) are also satisfied, a \textbf{positive totally ordered monoid}.
\end{definition}

\section{Ordering Defined by `Divisibility' in Monoids With Cancelation}
Let $M$ be a (possibly) noncommutative monoid. We assume cancelation on the right and on the left
$$a+c=b+c\Rightarrow a=b,\quad c+a=c+b\Rightarrow a=b.$$
That would be, of course, a prerequisite for any possibility to \textit{embed $M$ into a group.}
\begin{definition}\label{def:d}
The \textbf{right divisibility preorder} (resp.\ \textbf{left divisibility preorder}) on $M$ is defined as
\begin{equation}
a\le_{\textrm{r}}b\,:\Leftrightarrow\exists\,c\,\,(b=a+c)
\end{equation}
(resp.
\begin{equation}
a\le_{\textrm{l}}b\,:\Leftrightarrow\exists\,c\,\,(b=c+a)
\end{equation}
\end{definition}
One clearly sees that the associativity of the operation implies these are indeed (positive) preorders -- satisfy (O) and (P), and each satisfies \textit{one} of the assertions of (C): for right divisibility $a\le_{\textrm{r}}b\Rightarrow c'+a\le_{\textrm{r}}c'+b$, while for left divisibility $a\le_{\textrm{l}}b\Rightarrow a+c'\le_{\textrm{l}}b+c'$.
\begin{proposition}
The following are equivalent:
\begin{enumerate}
\item  \label{i:1} $M$ with right divisibility is a partial ordering, i.e.\ satisfies (A).
\item  \label{i:2} $M$ with left divisibility is a partial ordering, i.e.\ satisfies (A).
\item  \label{i:3} If a sum $a_1+\ldots+a_n$ equals $0$ then all the summands are $0$.
\end{enumerate}
\end{proposition}
\begin{Prf}
(\ref{i:1}) implies (\ref{i:2}) (and similarly (\ref{i:1}) implies (\ref{i:3})). Indeed, in (\ref{i:3})
\begin{eqnarray*}
&&0=a_1+\ldots+a_n\ge_{\textrm{r}}a_1+\ldots+a_{n-1}\ge_{\textrm{r}}\cdots\ge_{\textrm{r}}a_1\ge_{\textrm{r}}0\\
&&0=a_1+\ldots+a_n\ge_{\textrm{l}}a_2+\ldots+a_n\ge_{\textrm{l}}\cdots\ge_{\textrm{l}}a_n\ge_{\textrm{l}}0.\\
\end{eqnarray*}
On the other hand, assume (\ref{i:3}) holds.

Now, $a\le_{\textrm{r}}b\le_{\textrm{r}}a$ means $\exists c,c'$ such that $b=a+c$,\,\,$a=b+c'$. Then $a+0=a=a+c+c'$ and by cancelation $c+c'=0$, hence by (\ref{i:3}) $c=c'=0$, thus $a=b$. Similarly for $\le_{\textrm{l}}$.
\end{Prf}

\subsection{Extending the Monoid by `Annexing' Differences (Assuming the Ordering $\le_{\textrm{r}}$ Total)}
By cancelation the $c$ in Def.\ (\ref{def:d}), when it exists, is unique.

Thus if $a\le_{\textrm{r}}b$, we write the unique $c$ satisfying $b=a+c$ as $-a+b$.

And if $a\le_{\textrm{l}}b$, we write the unique $c$ satisfying $b=c+a$ as $b-a$.

(Note that $-(\,)$ is used here as $(\,)^{-1}$ is used for multiplicative notation -- minding the order of summands. In particular $-(a+b)=(-b)+(-a)$!)\pa

And in analogy with the way one formally extends the natural numbers to the integers or rationals, we may try to `annex' differences $b-a$. Yet with noncommutativity complications may well arise. (Compare  \cite{H} Ch.\ 7 about
quotient noncommutative rings and the Ore condition.) \pa

To fix matters, suppose \textbf{the ordering $\le_{\textrm{r}}$ is total}.\pa

Then define the \textbf{extension by differences of $M$} as the set $M\times M$, a pair $(a,b)$ there written as $b-a$, modulo the equivalence relation (stemming from an identity $(b+c)-(a+c)=b-a$ valid in some embedding group, say):
\begin{equation}
b_1-a_1\,\equiv\,b_2-a_2\,\Leftrightarrow\,\exists\,c,c'\,(b_1+c=b_2+c',\,\,a_1+c=a_2+c').
\end{equation}
That factor set (where we denote by $[b-a]$ the equivalence class of $b-a$) is endowed by the operation defined by: \begin{equation}\label{eq:o}
[c-b]+[b-a]=[c-a].
\end{equation}
(\ref{eq:o}) makes sense. Indeed, to find $[b_1-a_1]+[b_2-a_2]$,

one either has $a_1\le_{\textrm{r}}b_2$, then there is $c=-a_1+b_2$,\,\,$a_1+c=b_2$ and
$$[b_1-a_1]+[b_2-a_2]=[(b_1+c)-(a_1+c)]+[b_2-a_2]=[(b_1+c)-a_1]:$$
\begin{equation}
[b_1-a_1]+[b_2-a_2]=\left[\left(b_1+\left(-a_1+b_2\right)\right)-a_1\right]\quad\textrm{if }a_1\le_{\textrm{r}}b_2
\end{equation}
or $a_1\ge_{\textrm{r}}b_2$, then there is $c'=-b_2+a_1$,\,\,$a_1=b_2+c$ and
$$[b_1-a_1]+[b_2-a_2]=[b_1-a_1]+[(b_2+c)-(a_2+c)]=[b_1-(a_2+c)]:$$
\begin{equation}
[b_1-a_1]+[b_2-a_2]=\left[b_1-\left(a_2+\left(-b_2+a_1\right)\right)\right]\quad\textrm{if }a_1\ge_{\textrm{r}}b_2.
\end{equation}
By (\ref{eq:o}) this operation is associative. Also $[0-0]$ is a unit element and $[a-b]$ is the inverse of $[b-a]$, indeed $[b-a]+[a-b]=[b-b]=[0-0]$. The extension by differences is \textit{a group} $G$.\pa

One may identify an element $c\in M$ with $[c-0]$, the operation carrying over. indeed
$[c-0]+[c'-0]=[(c+c')-(0+c')]+[c'-0]=[(c+c')-0]$. $M$ is a submonoid of the group $G$.\pa

And conversely, given a submonoid $M$ of a group $G$, for the right divisibility preorder $M$ to be a total ordering one obviously needs
\begin{equation}
M\cap(-M)=\{0\},\qquad (-M)+M\subset M\cap(-M).
\end{equation}

\subsection{A Really Noncommutative Example}
Start with a kind of \textit{Heisenberg group $H$ over the integers}. Namely, $H$, as a set, is $\mathbb{Z}^3$, its elements written as $(m,n,p)$,\,\,$m.n.p\in\mathbb{Z}$,
 equivalently as $mi+nj+pk$ ($i,j,k$ being a basis). But in the
operation one introduces a `twist': $k$ commutes with everything, but while $i+j$ is as usual, $j+i$ is defined as
$i+j+k$, that making $mi+nj$ be as usual but $nj+mi:=mi+nj+(mn)k$.\pa

In $H$ we take a lexicographic ordering (and $M$ will be the set of elements $\ge0=(0,0,0)$), as follows:
\begin{equation}\label{eq:lx}
mi+nj+pk\le m'i+n'j+p'k:\Leftrightarrow m<m'\textrm{ or }(m=m'\textrm{ and }n<n')
\textrm{ or }(m=m'\textrm{ and }n=n'\textrm{ and }p\le p')
\end{equation}
So $M$, the set of elements $\ge0=(0,0,0)$ will be
\begin{equation}\label{eq:tw}
M=\left\{mi+nj+pk\,\big|\,m>0\textrm{ or }(m=0\textrm{ and }n>0)
\textrm{ or }(m=n=0\textrm{ and }p\ge0)\right\}
\end{equation}
And it is clearly stable with respect to $+$ -- the `twist' may affect only $p$, the coefficient of $k$, and that only if there is an $i$ term, but then $m>0$ and $p$ does not matter as per (\ref{eq:tw}).\pa

$M$ is not commutative: $i$ and $j$ belong to $M$ and $j+i=i+j+k\ne i+j$. We claim that still the right divisibility ordering there is a total ordering. Indeed,
\begin{proposition}
For $mi+nj+pk,\,m'i+n'j+p'k\in M$,\,\,$mi+nj+pk\le_{\textrm{r}}m'i+n'j+p'k$ (right divisibility ordering)
if and only if $mi+nj+pk\le m'i+n'j+p'k$ in the lexicographic ordering (\ref{eq:lx}).
\end{proposition}
\begin{Prf}
For $a=mi+nj+pk\le_{\textrm{r}}b=m'i+n'j+p'k$, there must be an $m''i+n''j+p''k\in M$, (i.e.\ either $m''>0$ or $m''=0,\,n''>0$ or $m''=n''=0,p''\ge0$) such that
\begin{eqnarray*}
&&b=m'i+n'j+p'k=a+m''i+n''j+p''k=mi+nj+pk+m''i+n''j+p''k\\
&&=(m+m'')i+(n+n'')j+(p+p''+nm'')k.
\end{eqnarray*}
If $m<m'$ or $m=m'$ and $n<n'$, i.e.\ the terms with $i$ and $j$ are not the same in $a$ and $b$, then the term with $k$ is of no concern and clearly the assertion holds. If these are the same, then $m''=n''=0$ and the `twist' vanishes -- $a$ and $b$ have $k$ terms $pk$ and $(p+p'')k$ and $p''\ge0$ and again the assertion holds.
\end{Prf}

\section{Ordering Compatible on the Right and on the Left}
From now on we insist on compatibility --  (C) satisfied -- on the right and on the left. Cancelation is not postulated.

\subsection{Using The Free Noncommutative Monoid}
Recall, that the \textbf{free noncommutative monoid} $\mathcal{F}_S$ over a set $S$ (referred to as \textit{alphabet} and its members as \textit{letters}) is the set of words over the alphabet $S$ with the operation of concatenation (including the empty word which serves as the unit element). Here it is written multiplicatively.
\pa
$\mathcal{F}_S$ satisfies the usual \textit{universal property}: Denote by $i:S\to\mathcal{F}_S$ the mapping sending each letter $s\in S$ to the one-letter word $s$. We say that by $i$,\,\,`$S$ is a subset of' $\mathcal{F}_S$. Then for any monoid $M$ and any mapping $f:S\to M$, $f$ has a unique `extension' $\mathcal{F}_f:\mathcal{F}_S\to M$, i.e.\ such that $\mathcal{F}_f\circ i=f$.
\pa
Now endow $\mathcal{F}_S$ with the following ordering $\preceq$: For words $v,w\in\mathcal{F}_S$, $v\preceq w$ if $v$ obtains from $w$ by deleting some of the letters, keeping the order of the letters not deleted. Call $\mathcal{F}_S$
with this ordering the \textbf{free positive (partially) ordered monoid} over $S$.
\pa
Then one easily sees that (O), (P), (C) and (A) are satisfied -- $\mathcal{F}_S$ becomes a positive partially ordered monoid, but in general not total -- (T) is not satisfied. Also if $M$ is a positive preordered monoid, then for any $f$ as above $\mathcal{F}_f$ will be order-preserving.
\pa
Therefore, assuming $M$ just a monoid and $\le$ a relation making it into a positive preordered monoid, then for the alphabet $M$ and the identity map $\textrm{Id}:M\to M$, (which, by the way, would imply the analogous facts for any alphabet $S$ and $f:S\to M$), \textit{$\le$ will be an extension of the push by $\mathcal{F}_\textrm{Id}$ of the ordering in $\mathcal{F}_M$}, i.e.\ of the relation
$$\mathcal{R}:=\left\{a,b\in M\times M\,\Big|\,\exists v,w\in \mathcal{F}_M\,\left(v\preceq w, a=\mathcal{F}_\textrm{Id}(v), b=\mathcal{F}_\textrm{Id}(w)\right)\right\},$$
hence also will be an extension of the transitive closure of $\mathcal{R}$, i.e.,
$$\left\{a,b\in M\times M\,\Big|\,\exists n \textrm{ and }c_0,c_1\ldots,c_n\left(a=c_0, c_0\mathcal{R}c_1, c_1\mathcal{R}c_2,\ldots,c_{n-1}\mathcal{R}c_n, b=c_n\right)\right\}.$$
That transitive closure, a relation defined canonically in any monoid $M$, which we denote by $\preceq_{\min}$, always
satisfies (O) of course, and also, as easily seen, (P) and (C), hence makes $M$ into a positive preordered monoid. And we saw that any relation so making $M$ is \textit{an extension of $\preceq_{\min}$}. But if the latter satisfies (A) -- is antisymmetric then so must be $\preceq_{\min}$. Thus,
\begin{proposition}\label{prop:1}
If in a monoid $M$ there exists any relation making it into a positive partially ordered monoid, then also $\preceq_{\min}$ in $M$ must satisfy (A).

And if so, $M$ with $\preceq_{\min}$ is a positive partially ordered monoid.
\qed
\end{proposition}
And $M$ can become a \textit{positive totally ordered monoid} only by some extension of $\preceq_{\min}$ that would satisfy (O), (C), (T), (A) ( (P) always `inherited' from $\preceq_{\min}$).
\pa
\begin{Rmk}
If the monoid $M$ we started with was \textit{commutative}, then one easily finds that the relation $\mathcal{R}$ is just $\exists c \left(b=a+c\right)$ which is transitive hence identical with its transitive closure $\preceq_{\min}$. That is, of course, the usual way to try to define a positive ordering in a commutative monoid.
\end{Rmk}
\subsection{The Axiom Partially Ordered Monoids Must Satisfy}
Now, to be more explicit about $\preceq_{\min}$, define:
\begin{definition}
Let $M$ be a (possibly noncommutative) monoid. A finite sequence of elements in $M$ will be called a \textbf{vector} (\textrm{being just a word in $\mathcal{F}_M$ as above}), with the number of elements its \textbf{length} and the sum of its elements (in the given order) its \textbf{weight} (\textrm{being just the above $\mathcal{F}_\textrm{Id}$ applied to it}). Two vectors with the same weight will be called \textbf{isobaric}. A finite sequence of vectors (in general of different lengths) will be called a \textbf{table}, with the vectors its \textbf{rows}.
\pa
An \textbf{augmentation} of a vector $v=(a_1,a_2,\ldots,a_n)$ is a vector obtained from $v$ by inserting some added element $c$: $v'=(a_1,a_2,\ldots,a_k,c,a_{k+1},\ldots,a_n)$,\,\,$k=0,1,\ldots,n$.
\pa
A table $A$ will be called \textbf{monotone} if each row in $A$ is isobaric with an augmentation of the previous row.
\end{definition}
Now, as one easily finds,
\begin{proposition}\label{prop:le1}
For elements $a,b$ in a (possibly noncommutative) monoid $M$, $a\preceq_{\min}b$ if and only if there is a monotone table $A$, such that the weight of its first row (resp.\ the last row) is $a$ (resp. $b$).
\qed
\end{proposition}
As for the requirement in Prop.\ \ref{prop:1} that $\preceq_{\min}$ in $M$ satisfy (A), it says that
$$a_1\mathcal{R}a_2,a_2\mathcal{R}a_3,\ldots,a_{n-1}\mathcal{R}a_n,a_{n}\mathcal{R}a_1\Rightarrow a_1=a_2=\cdots=a_n,$$
that is, that the following axiom holds:
\begin{axiom}\label{ax:le}
If a table $A$ is monotone and its the first and last row are isobaric -- have the same weight, then the weights of all the rows of $A$ are the same.
\end{axiom}
If this axiom holds for a (possibly noncommutative) monoid $M$ we shall say that the $M$ is (positively) \textbf{orderable}.
As wee saw, then $M$ with $\preceq_{\min}$ is a positive partially ordered monoid.
\pa
From now on we always assume our monoids are orderable unless specified otherwise. Of course, as we have seen,
when we speak about a \textit{partially ordered monoid}, in particular a \textit{totally ordered} one, it is automatically orderable. i.e.\ Axiom \ref{ax:le} holds.

\subsection{Multiplying by a Nonnegative Integer}
Let $M$ be a (positive) partially ordered (possibly noncommutative) monoid.

In $M$, as in any monoid, one can multiply a natural number $n$ by an $a\in M$:
\begin{definition}\label{def:1}
For $a\in M$ and $n=1,2,\ldots$ define:

$na:=a+\ldots+a$ ($n$ times $a$). Define also $0a=0$.
\end{definition}
By associativity the definition is unequivocal, and we have
\begin{equation}\label{eq:n}
(n+m)a=na+ma,\,\,(mn)a=m(na).
\end{equation}
Also, by Def.\ \ref{def:1}\,\,$m\le n\Rightarrow ma\le na$, and by (P) and (C) $a\le b\Rightarrow na\le nb$.\pa

\begin{corollary} (to Axiom \ref{ax:le})\label{cor:1}
For $a\in M$ and $n_1<n_2<n_3$, if $n_1a=n_3a$ then they are equal also to $n_2a$.
\end{corollary}
\begin{Prf}
Take in Axiom \ref{ax:le} a table $A$ with $n_3-n_2+1$ rows all whose entries are $a$, the $i$'th row of length
$n_1+i-1$. $A$ is obviously monotone. Since $n_1a=n_3a$ its first and last row are isobaric, thus by Axiom \ref{ax:le} the weights of all the rows are the same, hence the assertion.
\end{Prf}
\subsection{Absorbing Elements  and Idempotents}
\begin{definition}
$a$ \textbf{absorbs} $b$ (and $b$ is absorbed by $a$) \textbf{on the left} if $a+b=a$.

Similarly, a \textbf{absorbs} $b$ (and $b$ is absorbed by $a$) \textbf{on the right} if $b+a=a$.

$a$ is an \textbf{idempotent} if $a+a=a$ (that is, $a$ absorbs itself on the left/right).

$a$ is a \textbf{generalized idempotent} if $na=ma$ for some natural numbers $n\ne m$.
\end{definition}
Thus, if $a$ is not a generalized idempotent, than $n<m\,\,\Rightarrow\,\,na<ma$.\pa

Clearly, if $pa$ is a generalized idempotent for some $p=1,2,\ldots$ then so is $a$.\pa

\begin{Rmk}\label{rm:gi}
Suppose $a$ is a generalized idempotent. Then $\exists n<m\,\,na=ma$ and we take the pair $n,m$ with the \textit{least possible $n$}. This means that $a,2a,\ldots,na$ are all different.\pa

Then by Cor.\ \ref{cor:1} $na,(n+1)a,\ldots,ma$ are all equal, and by (\ref{eq:n}) also $(n+k)a,\dots,(m+k)a$ are equal for $k=1,2,\ldots$. And since for $k$ and $k+1$ these sequences overlap, we finally find that
$$na=(n+1)a=(n+2)a=\ldots.$$
\end{Rmk}
\begin{definition}
Let $P$ be a partially ordered set.

A subset $S\subset P$ is called \textbf{lower},\,\,(resp.\ \textbf{upper}) \textbf{hereditary} if $a\in S,\,b\le a\Rightarrow b\in S$ (resp.\ $a\in S,\,b\ge a\Rightarrow b\in S$).

For $a\in S$, the \textbf{lower cone},\,\,(resp.\ \textbf{upper cone}) with vertex $a$ is the set $\left\{b\in S\,\big|\,b\le a\right\}$,\,\,(resp.\ $\left\{b\in S\,\big|\,b\ge a\right\}$).
\end{definition}
In particular, in a (positive) partially ordered (possibly noncommutative) monoid $M$ we may speak of \textit{lower hereditary submonoids}. Then for an element $a\in M$, the lower hereditary submonoid generated by $a$ is easily seen to be
\begin{equation}\label{eq:lhm}
\left\{b\in M\,\big|\,\exists n=1,2,\ldots\,\left(b\le na\right)\right\}.
\end{equation}
Also, the sets of elements absorbed on the left (resp.\ on the right) by a fixed element $a$ in $M$,
\begin{eqnarray*}
&&\textrm{lAb}(a):=\left\{b\in M\,\big|\,a+b=a\right\},\\
&&\textrm{rAb}(a):=\left\{b\in M\,\big|\,b+a=a\right\},
\end{eqnarray*}
are clearly submonoids and by properties (P) (C) and (A) of the ordering also lower hereditary.\pa

Since one always has $b\le a+b,\,b+a$, we conclude that if $b$ is absorbed on the left or on the right by $a$ then $b\le a$.\pa

So $\textrm{lAb}(a)$ and $\textrm{rAb}(a)$ are lower hereditary sumonoids which are contained in the lower cone with vertex $a$ -- the latter, of course, lower hereditary but need not be a submonoid (i.e.\ need not contain the sum of two of its elements).\pa

And we clearly have
\begin{proposition}\label{prop:ab}
Let $M$ be a (positive) partially ordered (possibly noncommutative) monoid.

Then for some fixed $a\in M$ the following are equivalent:
\begin{enumerate}
\item $\textrm{lAb}(a)$ is the whole lower cone with vertex $a$, that is, $a$ absorbs on the left all the elements $\le$ from it.
\item $\textrm{rAb}(a)$ is the whole lower cone with vertex $a$, that is, $a$ absorbs on the right all the elements $\le$ from it.
\item $a$ is an idempotent: $a+a=a$, i.e.\ $a$ belongs to $\textrm{lAb}(a)$, equivalently belongs to $\textrm{rAb}(a)$.
\item The lower cone with vertex $a$ is a (lower hereditary) submonoid.
\item $\textrm{lAb}(a)$ has a greatest element, i.e.\ $\ge$ than all its other elements (necessarily equal to $a$).
\item $\textrm{rAb}(a)$ has a greatest element, i.e.\ $\ge$ than all its other elements (necessarily equal to $a$).
\end{enumerate}
\qed
\end{proposition}
\begin{Rmk}
To see that there are positive partially ordered (possibly) noncommutative monoids where some implications will not hold, say where for some $a$ and $b$,\,\,\,$a+b$ is a generalized idempotent while $a$ and $b$ are not, consider the following construction:\pa

Start with the free noncommutative momoid $\mathcal{F}_S$ over an alphabet $S$, partially ordered by $\preceq$ as above. Let $C$ be some upper hereditary subset of $\mathcal{F}_S$, say the upper cone of some $w\in\mathcal{F}_S$. Now take the equivalence relation $\approx$ in $\mathcal{F}_S$ {\em just collapsing $C$ to a point $[C]$}. Note that then $\mathcal{F}_S/{\approx}$ \textit{will always be a (positive partially ordered) monoid} -- the monoid operation there well-defined, making $[C]$ a greatest element there (hence absorbing all others on the left and on the right) and \textit{with no relations between words on $S$ except the equality among all members of $C$}.\pa

For example: if $a,b\in S$ and $C$ is the upper cone of the word $(a,b)$, then $a+b$ is an idempotent in $\mathcal{F}_S/{\approx}$ while neither of $a$, $b$ is a generalized idempotent!\pa

On the other hand one notes that since $n(b+a)=b+(n-1)(a+b)+a$ we do have always: {\em for $a$ and $b$ in a monoid $M$, $b+a$ is a generalized idempotent if and only if $a+b$ is so}.\pa

One may loosely say that when multiplying by a large natural number, {\em the difference between $a+b$ and $b+a$ somewhat dwindles}.  We shall pursue this clue in \S\ref{sc:ra}.
\end{Rmk}

\section{The Case of Total Ordering, Ratios}
Let $M$ be a (positive) totally ordered (possibly noncommutative) monoid.\pa

Then a lower (resp.\ upper) hereditary subset is a \textit{lower ray} (\textit{upper ray}) (possibly empty or the whole $M$) and a set is a lower ray if and only if its complement is an upper ray.

\subsection{Archimedicity and Possible Absorbing}
\begin{definition}
For $a,b\in M$, $a$ is said to be \textbf{infinitely greater than $b$}, (and $b$ \textbf{infinitely smaller than $b$}) if there is no $n=1,2,\ldots$ such that $a\le nb$, in other words, $a>nb$ for all $n$. Otherwise, i.e.\ if there exists an $n$ with $a\le nb$, then $a$ is called \textbf{Archimedean} to $b$.

If $a$ and $b$ are both Archimedean to each other, we say that $a$ and $b$ are \textbf{commensurable}.
\end{definition}
\begin{Rmk}
Clearly the following relations are transitive: being infinitely greater than; being Archimedean to; being commensurable to. Also, for a fixed $a\in M$,\,\,$M$ is partitioned into three sets (possibly empty): the lower ray of the $b$'s infinitely smaller then $a$;, all these smaller then the $b$'s commensurable to $a$; -- these two sets together make those to which $a$ is Archimedean;, and all these smaller than the $b$'s infinitely greater than $a$ which make an upper ray.
\end{Rmk}
And, (cf.\ the previous \S), for an $a\in M$, the lower hereditary submonoid generated by $a$ is
$$\left\{b\in M\,\big|\,\exists n=1,2,\ldots\,\left(b\le na\right)\right\},$$
thus it coincides with the set of elements of $M$ which are Archimedean to $a$.\pa

In particular one notes that the set of elements Archimedean to $a$ is a (lower hereditary) submonoid -- contains the sum of any two of its members.\pa

Then we have:
\begin{proposition}
For  $a,b\in M$,\,\,$a$ absorbing $b$ on the left or right is possible only if either

$a$ is an idempotent: $a=a+a$ (hence absorbs all the elements $\le$ from it), or

$a$ is infinitely greater than $b$.
\end{proposition}
\begin{Prf}
By the previous \S, the set of elements absorbed by $a$ on the left (resp.\ on the right) (which is contained in the lower cone with vertex $a$) is a lower hereditary monoid, hence contains the lower hereditary monoid generated by $b$ which is the set of elements Archimedean to $b$.

Then either the latter is the whole lower cone with vertex $a$, equivalently contains $a$: $a$ is an idempotent.

Or it does not contain $a$, thus $a$ is not Archimedean to $b$: $a$ is infinitely greater than $b$.
\end{Prf}
So
\begin{proposition}\label{prop:ar}
For commensurable $a,b\in M$,\,\, $a$ absorbing $b$ on the left or right is possible only if
$a=nb$ for some natural $n$, so we are in the situation of Remark \ref{rm:gi}, in particular $b$ is a generalized idempotent.
\end{proposition}
\begin{Prf}
By the previous proposition, $a$ must be an idempotent. Also since $a$ is also Archimedean to $b$, there is a natural number $n$ such that $a\le nb$. Let $n$ be the smallest satisfying that. If $n=1$ than $b\le a\le b\Rightarrow a=b$. If $n>1$ then $b\le a$, so $a$ absorbs $b$ hence also absorbs $nb$ thus one has $a\ge nb$ and finally $a=nb$.
\end{Prf}
And moreover
\begin{proposition}
If $a$ and $b$ are commensurable, then if one of them is a generalized idempotent so is the other, and in fact there are natural numbers $n$ and $m$ such that $a,2a,\ldots,na$ are different, $b,2b,\ldots,nb$ are different, while $na$ is equal to $mb$ and $na=(n+1)a=(n+2)a=\ldots=mb=(m+1)b=(m+2)b=\ldots$.

And the latter element, an idempotent, is the greatest in the lower ray of all elements to which $a$ (equivalently $b$) is Archimedean, and absorbs them on the right and on the left.
\end{proposition}
\begin{Prf}
Suppose $a$ is a generalized idempotent and $n$ as in Remark \ref{rm:gi}. Then the `maximal' $na$ is an idempotent, absorbing itself hence elements smaller than it on the left and on the right.

Now $b$ is commensurable with $a$. Therefore there is some $n'$ so that $b\le n'a$, thus $b$, hence all its multiples, are less or equal than the maximal idempotent $na$. On the other hand $a$, hence $na$, is commensurable with $b$. So there is an $m$ such that $na\le mb$, take the minimal with this property. But then $mb=na$ and the assertion follows.
\end{Prf}

\subsection{Ratios \`{a} la Euclid's \textit{Elements}}\label{sc:ra}
Now let us try to follow, in our setting, what Euclid's \textit{Elements} does. There, in order to impose numerical values on geometrical objects which, at the start, one just adds and compares, a theory of `ratios' (due to Eudoxos) is developed, a most abstract and almost `modern' by the standards of the \textit{Elements}. (Where they already knew that they cannot expect all ratios to be `rational' -- given by a fraction of two integers -- the Pythagoreans had discovered, two centuries before, that the ratio of the diagonal of a square to its side is irrational!)\pa

So, let $a,b\in M$, commensurable and not generalized idempotents. Then by Prop.\ \ref{prop:ar} none of them can absorb the other.\pa

Let $m,n$ be positive integers. We put a sign $=,\ge,\le,>,<$ between \textit{the ratio $a:b$} and \textit{the ratio $m:n$} if that sign holds between $na$ and $mb$.\pa

Now, if $a:b=m:n$, i.e., $na=mb$ then we have for any positive integer $p$,\,\,$(pn)a=p(na)=p(mb)=(pm)b$, meaning $a:b=(pm):(pn)$.\pa

And if $a:b\le m:n$, i.e., $na\le mb$, then $(pn)a=p(na)\le p(mb)=(pm)b$, i.e., $a:b\le(pm):(pn)$.\pa

This means that for a fraction $q=m/n$, the signs $a:b$ has with two ways to write it as a ratio between a numerator and a denominator, namely $(p_1m):(p_1n)$ and $(p_2m):(p_2n)$, may certainly not be $<$ vs.\ $>$ or $>$ vs. $<$ -- they are both $\le$ or $=$ (write then $a:b\le q$) or both $\ge$ or $=$ (writing $a:b\ge q$).\pa

Now suppose the fractions $q'<q$. Write them with a common denominator as $q'=m'/n,\,q=m/n,\,\,m'<m$. Then if
$a:b\le q'$, i.e., $na\le m'b$, then, since $m'b<mb$, ($b$ is not a generalized idempotent!), we have  $na<mb$, i.e., $a:b<m:n$, so \textit{$a:b\ge q$ does not hold}. Similarly if $a:b\ge q$ then we find $a:b>m':n$ hence $a:b\le q'$ does not hold.\pa

So we conclude that, given $a,b\in M$, we are in the situation featuring with the \textit{Elements}' `ratios'. Namely, the set of all positive rational numbers $q$ is divided into two sets: those where $a:b\ge q$ and those where $a:b\le q$, and these form a \textit{Dedekind cut} -- both are nonempty since $a$ and $b$ are commensurable. The former is of the form $[\alpha,\to)$ or $(\alpha,\to)$ and the latter is of the form $(0,\alpha]$ or $(0,\alpha)$, for some real $\alpha>0$ which `the cut determined'.\pa

And we define this $\alpha$ as \textbf{the ratio $a:b$}. Thus $a:b\in(0,+\infty)$.

\subsection{Properties of Ratios}
\begin{proposition}
For $a,b,c\in M$ commensurable, not generalized idempotents,
\begin{enumerate}
\item $a:a=1$.\label{i:un}
\item $(b+a):(a+b)=1$\label{i:com}
\item $(na):a=n$\label{i:n}
\item $a:c=(a:b)\cdot(b:c)$.\label{i:tr}
\item $b:a=(a:b)^{-1}$.\label{i:inv}
\item $a\le b\Rightarrow a:c\le b:c$ and $c:a\ge c:b$.\label{i:le}
\item $(a+b):c=(b+a):c=(a:c)+(b:c)$.\label{i:ad}
\item If $a$ and $b$ are commensurable and not generalized idempotents, but $c$ is infinitely smaller than $a$ (equivalently than $b$), then $(c+a):b=(a+c):b=a:b$,\quad$a:(c+b)=a:(b+c)=a:b$.
\item In particular, also in the case that $c$ is infinitely smaller than $a$ we have $(a+c):(c+a)=1$.
\item Conversely, for $b$ Archimedially related to $a$, thus not infinitely smaller than it, such equalities cannot occur. Indeed then $(a+c):a>1$,\,$(c+a):a>1$,\,$a:(a+c)<1$ and $a:(c+a)<1$.
\end{enumerate}
\end{proposition}
\begin{Prf}
\begin{enumerate}
\item Obviously, if $n\le m$ then $na\le ma$. Hence $n:m\le a:a$ and $m:n\ge a:a$. Therefore $n/m>a:a$ cannot hold and $m/n<a:a$ cannot hold. Having that whenever $n\le m$, $a:a$ must be $1$.
\item We have $(n+1)(a+b)=a+n(b+a)+b\ge n(b+a)$. Thus $(n+1):n\ge(b+a):(a+b)$ and similarly $(n+1):n\ge(a+b):(b+a)$, i.e.\ $n:(n+1)\le(b+a):(a+b)$. Therefore $(b+a):(a+b)$ cannot be $>(n+1)/n$ and cannot be $<n/(n+1)$. That being the case for every $n$, $(b+a):(a+b)$ must be $1$.
\item Suppose $(na):c<n$ (resp.\ $(na):a>n$). Then there must be natural numbers $m,p$ such that $(na):a>p/m,\,p/m>n$ (resp.\ $(na):a<p/m,\,p/m<n$). This means that $p>mn$ while $(mn)a\le pa$ is impossible (resp.\ $p<mn$ while $(mn)a\ge pa$ is impossible). But, $p>mn$ surely implies $(mn)a\le pa$ (resp.\ $p<mn$ surely implies $(mn)a\ge pa$) -- contradiction.
\item Suppose $a:c<(a:b)\cdot(b:c)$ (resp.\ $a:c<(a:b)\cdot(b:c)$). Then there must be natural numbers $n,m,p$ such that $a:c>p/n,\,a:b<m/n,\,b:c<p/m$ (resp.\ $a:c<p/n,\,a:b>m/n,\,b:c>p/m$). This means that $na\le pc$ is impossible and so are $na\ge mb$ and $mb\ge pc$ (resp.\ $na\ge pc$ is impossible and so are $na\le mb$ and $mb\le pc$). But, the ordering being total, that means $na>pc>mb>na$ (resp.\ $na<pc<mb<na$) -- impossible.
\item The assertion follows, of course, from items \ref{i:tr} and \ref{i:un}, but one may give a direct proof:

     It is enough to show that for any $n,m$, $n/m$ less than (resp.\ greater than) $b:a$ implies the same (even if in the weak sense) with $(a:b)^{-1}$. But
\begin{eqnarray*}
&& n/m<b:a \Rightarrow n:m\textrm{ cannot be}\ge b:a \Rightarrow na\ge mb\textrm{ does not hold }\\
&& \Rightarrow m:n\le a:b\textrm{ does not hold }\Rightarrow m/n \ge a:b \Rightarrow n/m \le (a:b)^{-1}.
\end{eqnarray*}
    And similarly for the reverse inequalities.
\item The second assertion follows from the first and item \ref{i:inv}. For the first assertion, we have to show that $n/m$ less than $a:c$ implies the same (even if in the weak sense) with $b:c$. Indeed
\begin{eqnarray*}
&& n/m<a:c \Rightarrow n:m\textrm{ cannot be}\ge a:c \Rightarrow nc\ge ma\textrm{ does not hold }\\
&& \Rightarrow\textit{ a fortiori}\textrm{, as }a\le b
\textrm{ hence }ma\le mb,\,\,\,nc\ge mb\textrm{ does not hold }\\
&&\Rightarrow n:m\ge b:c\textrm{ does not hold }\Rightarrow n/m \le b:c.
\end{eqnarray*}
\item Here, in order to partially circumvent the noncommutativity, we rely on the results of the previous items. Take $c$ fixed and consider the map $\rho_c:a\mapsto a:c$ between the set of elements in $M$ commensurable to $c$ and $(0,\infty)\subset\mathbb{R}$, the latter endowed with addition (of numbers). We wish to prove that it preserves the operation.

    By the previous items we know that, for variable $a,b\in M$,\,\,$\rho_c(c)=1$,\,\,
    $a:b=\rho_c(a)/\rho_c(b)$,\,\,and $a\le b\Rightarrow\rho_c(a)\le\rho_c(b)$ -- $\rho_c$ is monotone. Also by item \ref{i:n}\,\,$\rho(na)=n\rho(a)$. And since $(b+a):(a+b)=1$ one obtains that $\rho_c(b+a)=\rho_c(a+b)$.

    Thus the value $\rho_c$ gives to a sum does not depend on the order of the summands. In particular $n\rho_c(a+b)=\rho_c(n(a+b))=\rho_c(na+nb)$

    So suppose that for some $a$ and $b$\,\,\,$(a+b):c<(a:c)+(b:c)$ (resp.\ $(a+b):c>(a:c)+(b:c)$).
    Then there must be natural numbers $n,m,p$ such that $a:c<m/n,\,b:c<p/n,\,(a+b):c>(m+p)/n$ (resp.\ $a:c>m/n,\,b:c>p/n,\,(a+b):c<(m+p)/n$). This means that $n(a+b)\le(m+p)c$ is impossible and so are $na\ge mc$ and $nb\ge pc$ (resp.\ $n(a+b)\ge(m+p)c$ is impossible and so are $na\le mc$ and $nb\le pc$). The ordering being total, that means $na+nb\le(m+p)c<n(a+b)$ (resp.\ $na+nb\ge(m+p)c>n(a+b)$). Applying $\rho_c$ which is monotone and gives to both $na+nb$ and $n(a+b)$ the same value $n\rho_c(a+b)$ and to $(m+p)c$ the value $m+p$, we find $(a+b):c=\rho_c(a+b)=(m+p)/n$. That contradicting the strict inequality here that we had.
\item First, since $c$ is infinitely smaller than $a$, implying of course $c<a$, and also $a+c\ge c$, $a+c$ and $a$ are surely commensurable. By items \ref{i:tr} and \ref{i:inv} it suffices to prove our assertions for some fixed $b$.

    Thus let $a$ and $c$ vary and use the notation of item \ref{i:ad}, for some fixed $a_0$:
    $\rho_{a_0}(a):=a:a_0$ and also $\beta(a,c):=\rho_{a_0}(a+c)-\rho_{a_0}(a)$,\,\,$\beta'(a,c):=\rho_{a_0}(c+a)-\rho_{a_0}(a)$. Then, since $a+c,\,c+a\ge a$,\,\,$\beta(a,c),\,\beta'(a,c)\ge0$ and our task is to prove that they both vanish. We use the properties of $\rho$ described in \ref{i:ad}.

    First they yield that $\beta=\beta'$, both not depending on $a$. Indeed, for any $a$ and $b$,
        $$\beta(a,c):=\rho_{a_0}(a+c)-\rho_{a_0}(a)=\rho_{a_0}(a+c+b)-\rho_{a_0}(b)-\rho_{a_0}(a)=
        \rho_{a_0}(c+b)-\rho_{a_0}(b)=\beta'(b,c).$$
    Thus $\beta(a,c)=\beta'(a,c)=\beta(c)$.

    Now this $\beta$ is clearly \textit{additive in $c$}, since
    $$\beta(c+c')=\rho_{a_0}(a+c+c')-\rho_{a_0}(a)=\left(\rho_{a_0}(a+c+c')-\rho_{a_0}(a+c)\right)+
    \left(\rho_{a_0}(a+c)-\rho_{a_0}(a)\right)=\beta(c')+\beta(c).$$
    In particular $\beta(nc)=n\beta(c)$.

    But for any $a$ and $b$, commensurable while $c$ infinitely smaller than them, $c$ is surely \textit{smaller} than $b$. Therefore $a+c\le a+b$ making
    $$\beta(c)\le\rho_{a_0}(a+b)-\rho_{a_0}(a)=\rho_{a_0}(b),$$
    -- any value taken by $\rho_{a_0}$ at any $b$ bounding $\beta(c)$ for any $c$! But then it also bounds $\beta(nc)=n\beta(c)$ for any $n$. So $\beta$ must vanish.
\item Follows from the previous item.
\item That follows from item \ref{i:ad} and from $c:a$ being always \textit{positive}.
\end{enumerate}
\end{Prf}
So one may loosely say that \textit{for any (positive) totally ordered monoid, as far as commensurable elements, not generalized idempotents, are concerned, and these up to an `error' infinitely smaller than them, then the `ratio', as a `first approximation', imposes a numerical additive structure, thus commutative}.

\begin{acknowledgement}
I am much grateful to Professor Orr Shalit for posing and remarking on the problem which is the impetus for this note and on this note in preparation. Some ideas in this note are from him.
\end{acknowledgement}

\end{document}